\documentclass[12pt,a4paper]{amsart}
\usepackage{amssymb}
\usepackage{amsmath}
\usepackage{amsfonts}

\newcommand{\Q}{\mathbb{Q}}

\newcommand{\Z}{\mathbb{Z}}
\newcommand{\F}{\mathbb{F}}

\newtheorem{main-dummy}{Main-Dummy}

\newtheorem{main-theorem}[main-dummy]{Theorem}

\numberwithin{dummy}{section}
\numberwithin{equation}{section}

\newtheorem*{lemma*}{Lemma}

\newtheorem*{theorem*}{Theorem}

\newtheorem*{prop*}{Proposition}

\newtheorem*{cor*}{Corollary}

\theoremstyle{definition}

\theoremstyle{remark}

\newtheorem*{rem*}{Remark}

\begin{document}
\bibliographystyle{amsalpha}
\author{Sandro Mattarei}

\email{mattarei@science.unitn.it}

\urladdr{http://www-math.science.unitn.it/\~{ }mattarei/}

\address{Dipartimento di Matematica\\
  Universit\`a degli Studi di Trento\\
  via Sommarive 14\\
  I-38050 Povo (Trento)\\
  Italy}

\title{Exponential functions in prime characteristic}

\begin{abstract}
In this note we determine all power series $F(X)\in 1+X\F_p[[X]]$ such that
$(F(X+Y))^{-1} F(X)F(Y)$ has only terms of total degree a multiple of $p$.
Up to a scalar factor, they are all the series of the form
$F(X)=E_p(cX)\cdot G(X^p)$
for some $c\in\F_p$ and $G(X)\in 1+X\F_p[[X]]$,
where $E_p(X)=\exp\big(\sum_{i=0}^{\infty}X^{p^i}/p^i\big)$
is the Artin-Hasse exponential.
\end{abstract}

\date{30 June 2005}

\subjclass[2000]{Primary 39B50}

\keywords{functional equation, exponential series, Artin-Hasse exponential.}

\thanks{The  author  is grateful  to  Ministero dell'Istruzione, dell'Universit\`a  e
  della  Ricerca, Italy,  for  financial  support of the
  project ``Graded Lie algebras  and pro-$p$-groups of finite width''.}

\maketitle

\thispagestyle{empty}

\section{Introduction}

Even undergraduates probably know that the functional equation $f(x+y)=f(x)f(y)$ characterizes the exponential function $\exp(x)$,
up to a scaling in $x$, in large classes of real or complex functions.
For example, assuming $f$ strictly positive for simplicity, by taking logarithms of both sides one can reduce the equation
to Cauchy's equation $g(x+y)=g(x)+g(y)$, which has only homogeneous linear functions as solutions, under mild restrictions on $g$.
A slightly different route, if $f$ is differentiable, is to differentiate both sides of the equation with respect to $y$,
then evaluate for $y=0$, and finally solve the resulting differential equation $f'(x)=f(x)f'(0)$.
The last step is usually taught in calculus courses by obtaining a recurrence relation for the coefficients
of the Taylor expansion of $f$ at zero.
The latter proof applies also in a formal setting to show that, over any field $k$ of characteristic zero,
the formal exponential series $F(X)=\sum_{n=0}^{\infty}X^n/n!$ is the unique solution in $k[[X]]$ of the functional equation
\begin{equation}\label{eq:functional}\tag{$\ast$}
F(X+Y)=F(X)F(Y)
\end{equation}
satisfying $F'(0)=1$.

When $k$ is a field of positive characteristic, the series defining the exponential function does not make sense,
because most denominators vanish.
In fact, the argument sketched above proves that the only series $F(X)\in k[[X]]$ which satisfy~\eqref{eq:functional} in
$k[[X,Y]]$ are the constant series $0$ and $1$.
The absence of an analogue of the exponential function satisfying~\eqref{eq:functional} in characteristic $p$
is explained, in fancier terms, by the fact that the additive and multiplicative one-dimensional formal group laws
$\mathbf{\hat G}_\mathrm{a}(X,Y)=X+Y$ and $\mathbf{\hat G}_\mathrm{m}=X+Y+XY$ over $k$
are not isomorphic, see~\cite[(1.4.2)]{Haz}.
At this point we should mention that the {\em Carlitz exponential} (see~\cite{Hayes:Carlitz-Ch.19} for a pleasant introduction
to this topic) and its generalizations,
which are functions defined by power series with coefficients in function fields over finite fields, hence in positive characteristic,
display remarkable analogies with the complex exponential function;
however, those functions are additive, rather than multiplicative
as expressed by the functional equation~\eqref{eq:functional}.

Over $p$-adic fields a modified exponential function, the {\em Artin-Hasse exponential,}
defined by the series $E_p(X)=\exp\big(\sum_{i=0}^{\infty}X^{p^i}/p^i\big)$,
plays an important role.
This series turns out to have coefficients in the ring $\Z_p$ of $p$-adic integers,
and so has a larger radius of convergence than the ordinary exponential~\cite[(VII.2.2)]{Robert}.
In particular, $E_p(X)$ has a well-defined reduction modulo $p$, which we denote here by $\bar E_p(X)$.

Among many remarkable properties of the Artin-Hasse series there is one related with
the functional equation~\eqref{eq:functional}, namely, all terms of the series
$(E_p(X+Y))^{-1}E_p(X)E_p(Y)\in\Z_p[[X,Y]]$ have multiples of $p$ as total degrees.
The corresponding fact for the reduction $\bar E_p(X)$ of $E_p(X)$ modulo $p$
has already several applications.
It plays a role at some places in Gerstenhaber's theory of deformations of algebras,
in the case of a ground field of prime characteristic.
In particular, we mention~\cite[Theorem~2]{Ger:deformation1}, where a certain series related to $\bar E_p(X)$
(but $\bar E_p(X)$ would work as well)
is applied to a derivation $D\in Z^1(A,A)$ of an (associative) algebra $A$ to prove that its primary obstruction cocycle
$\mathrm{Sq}_p(D)\in Z^2(A,A)$
can be integrated to a (formal) deformation of $A$.
Another place is~\cite[pp.~59-60]{GerSch}, where, starting from a derivation $D$ of $A$ which is not integrable to order $p$,
the Artin-Hasse exponential is used to construct
a non-trivial (formal) deformation of $A$ which becomes trivial when the deformation parameter $t$ is replaced by $t^p$.
A more recent application of this property of $\bar E_p(X)$ appears in~\cite{Mat:Artin-Hasse},
where Artin-Hasse exponentials of derivations
are used for the construction of gradings of nonassociative algebras, in characteristic $p$, over cyclic groups of order $p^k$
with $k>1$.

The Artin-Hasse series can undergo substantial modifications without affecting the property considered in the previous paragraph.
(In fact, the series $\exp(X)\cdot\big(\sum_{m\ge 0} X^{mp}/(mp)!\big)^{-1}$
used by Gerstenhaber in~\cite[Theorem~2]{Ger:deformation1}
appears rather distant from the classical Artin-Hasse series.)
The goal of this note is to determine all such modifications.
More precisely, we prove the following result.

\begin{theorem*}
For a series $F(X)\in 1+X\F_p[[X]]$, the series
\[
(F(X+Y))^{-1}F(X)F(Y)\in\F_p[[X,Y]]
\]
has only terms of total degree a multiple of $p$ if, and only if,
\[
F(X)=\bar E_p(cX)\cdot G(X^p),
\]
for some $c\in\F_p$ and $G(X)\in 1+X\F_p[[X]]$.
\end{theorem*}

Note that the inverse $(F(X+Y))^{-1}$ makes sense only if $F(X)$ has a nonzero constant term,
and then it is not restrictive to assume that the constant term is $1$, as in the Theorem.
The sufficiency of the condition is almost obvious once we know that $\bar E_p(X)$ has the required property,
which results from a straightforward verification (as in the proof of~\cite[Theorem~2.3]{Mat:Artin-Hasse}, for example).

We present two slightly different proofs of the necessity.
One proof involves working over the ring $\Z_p$ of $p$-adic integers first,
as in the Proposition of the next section.
The other proof, more direct but perhaps less elegant, only involves series over $\F_p$.
Although the latter proof (and hence the Theorem) is valid over any field of positive characteristic in place of the prime field $\F_p$,
we only sketch it in a Remark and instead give the former proof in full.
We justify this choice on practical grounds.
On one hand, the case of series over the prime field is the one that matters for the applications mentioned earlier.
On the other hand, series satisfying the Theorem are usually explicitly produced as reductions modulo $p$
of series over the $p$-adic integers;
therefore, the explicit description which we give in the Proposition of all series in $\Z_p[[X]]$ whose
reduction modulo $p$ satisfies the condition of the Theorem might be of interest.

We conclude with another characterization of the series under study.
\begin{cor*}
For a series $F(X)\in 1+X\F_p[[X]]$, the series
\[
(F(X+Y))^{-1}F(X)F(Y)\in\F_p[[X,Y]]
\]
has only terms of total degree a multiple of $p$ if, and only if,
\[
F'(X)/F(X)=c\sum_{i=0}^\infty X^{p^i-1}.
\]
for some $c\in\F_p$.
\end{cor*}

\section{Proofs}

Before starting with the proofs of the Theorem and the Corollary
we note that the simplest variation on the Artin-Hasse series is replacing
it with any power series of the form
$\exp\big(\sum_{i=0}^{\infty}b_iX^{p^i}\big)$
having all coefficients in the ring $\Z_p$ of $p$-adic integers.
This is known to occur if and only if
$pb_i-b_{i-1}\in p\Z_p$ for all $i\ge 0$, where $b_{-1}=0$,
see~\cite[Proposition~4]{Die} or~\cite[Proposition~(2.3.3)]{Haz}.
The following result has a similar flavour.

\begin{prop*}\label{lemma:generalized-AH}
The following conditions on $F(X)\in 1+X\Q_p[[X]]$ are equivalent:
\begin{enumerate}
\item
$F(X)\in\Z_p[[X]]$, and all terms of the reduction modulo $p$ of
$(F(X+Y))^{-1}F(X)F(Y)$
have total degree a multiple of $p$;
\item
$F(X)=\exp\big(\sum_{j=1}^\infty c_jX^j\big)$, with
$pc_{pj}-c_j\in p\Z_p$ for all $j$, $c_1\in\Z_p$, and $c_j\in p\Z_p$ for $j>1$ and not a multiple of $p$.
\end{enumerate}
\end{prop*}

\begin{proof}
Since the series $F(X)$ has constant term $1$, its logarithm $G(X)=\log(F(X))\in\Q_p[[X]]$ is defined
by composition with the power series $\log(1+X)=\sum_{i=1}^\infty (-1)^{i+1}X^i/i$,
and we have
$F(X)=\exp(G(X))$.
Write
$G(X)=\sum_{j=1}^\infty c_jX^j$.

According to the Dieudonn\'{e}-Dwork criterion~\cite[(VII.2.3)]{Robert}, a power series $F(X)\in 1+X\Q_p[[X]]$
has coefficients in $\Z_p$ if and only if
\[
F(X)^p/F(X^p)\equiv 1\pmod{p},
\]
where the congruence sign means that the difference of the two sides belongs to $pX\Z_p[[X]]$.
Since $F(X)=\exp(G(X))$ here, it is convenient to state the criterion in terms of $G(X)$.
The condition of the criterion can be written as
\[
\exp\left(pG(X)-G(X^p)\right)\equiv 1\pmod{p}.
\]
Standard $p$-adic evaluations like in~\cite[(V.4.2)]{Robert} show
that composition on the left with $\exp$ and $\log$ gives inverse bijections between
$pX\Z_p[[X]]$ and $1+pX\Z_p[[X]]$.
Thus, we have the following exponential version of the Dieudonn\'{e}-Dwork criterion (see also~\cite[(VII.3.3)]{Robert}
or~\cite[2.2 and 10.2]{Haz} for far-reaching generalizations):
for a power series $G(X)\in X\Q_p[[X]]$, its exponential
$F(X)=\exp(G(X))$ belongs to $\Z_p[[X]]$ if and only if
\[
pG(X)\equiv G(X^p)\pmod{p}.
\]
In particular, the series $F(X)$ under consideration in this proof belongs to $\Z_p[[X]]$ if and only if
$pc_{pj}-c_j\in p\Z_p$ for all $j$, and $c_j\in\Z_p$ if $j$ is not a multiple of $p$.
Note that these conditions imply inductively that
$jc_j\in\Z_p$ for all $j$.

We introduce a further indeterminate $T$, and observe that the reduction modulo $p$ of a power series in $\Z_p[[T]]$
has only terms whose degrees are multiples of $p$
exactly when the derivative of the series belongs to $p\Z_p[[T]]$.
We apply this observation to $(F(TX+TY))^{-1}F(TX)F(TY)$, which belongs to
$\Z_p[[X,Y]]\,[[T]]$, but also to its subring $\Z_p[T]\,[[X,Y]]$, according to a natural identification.
Thus, for $F(X)\in\Z_p[[X]]$ the further condition in~(1) is equivalent to
\[
\frac{d}{dT}
\left((F(TX+TY))^{-1}F(TX)F(TY)\right)
\equiv 0\pmod{p},
\]
the congruence sign meaning now that the left-hand side belongs to $p\Z_p[T]\,[[X,Y]]$.
After performing the calculations and setting $T=1$, the condition can be rewritten as
\[
(X+Y)\frac{F'(X+Y)}{F(X+Y)}
\equiv
X\frac{F'(X)}{F(X)}+
Y\frac{F'(Y)}{F(Y)}
\pmod{p},
\]
that is,
\[
(X+Y)G'(X+Y)
\equiv
XG'(X)+
YG'(Y)
\pmod{p}.
\]

A power series $H(X)\in\F_p[[X]]$ is {\em additive,} that means, it satisfies
$H(X+Y)=H(X)+H(Y)$ in $\F_p[[X,Y]]$, exactly when all its terms have degrees which are powers of $p$.
This fact is well known at least for polynomials, see~\cite[Ch.~5, Sec.~2]{FesVos}, for example,
and it extends readily to power series upon noting that
$\F_p[[X]]/(X^{n})\cong\F_p[X]/(X^{n})$ for all $n$.
Hence our condition on $G(X)$ is equivalent to
\[
XG'(X)=\sum_{j=1}^{\infty}jc_jX^{j}
\in\Z_p[[X]]
\]
and $jc_j\in p\Z_p$ unless $j$ is a power of $p$.
That is,
$jc_j\in\Z_p$ if $j$ is a power of $p$, and
$jc_j\in p\Z_p$ otherwise.
It is now easy to see that these conditions, together with those found above ensuring that $F(X)\in\Z_p[[X]]$, amount to
$c_{pj}-c_j/p\in \Z_p$ for all $j$, $c_1\in\Z_p$, and $c_j\in p\Z_p$ for $j>1$ and not a multiple of $p$.
\end{proof}

\begin{proof}[Proof of the Theorem]
As we have mentioned in the Introduction, only the necessity of the condition requires a proof.
Let $F(X)\in 1+X\F_p[[X]]$, and suppose that
\[
(F(X+Y))^{-1}F(X)F(Y)\in\F_p[[X,Y]]
\]
has only terms of total degree a multiple of $p$.
Any series $\tilde F(X)\in 1+X\Z_p[[X]]$ which has $F(X)$ as its reduction modulo $p$
satisfies condition~(1) of the Proposition and, therefore, also condition~(2).
Since
\[
\exp\bigg(\sum_{p>1,\,p\nmid j} c_jX^j\bigg)\in\exp(pX\Z_p[[X]])=1+pX\Z_p[[X]],
\]
the reduction modulo $p$ of $\tilde F(X)$ is the same as that of the series
\[
\exp(c_1X)\exp\bigg(\sum_{j=1}^{\infty} c_{pj}X^{pj}\bigg)\in 1+X\Z_p[[X]].
\]
The quotient of the latter by $E_p(c_1X)$ belongs to $\Z_p[[X^p]]$.
Reduction modulo $p$ yields the desired conclusion.
\end{proof}

\begin{rem*}
We sketch a proof of the Theorem which works directly over $\F_p$, avoiding recourse to the $p$-adic integers.
Since the core of the argument is the same as in the proof of the Proposition, we will be very succinct.
Like in the proof of the Proposition we show that $XF'(X)/F(X)$ is an additive series, whence it has the form
$\sum_{i=0}^{\infty}b_iX^{p^i}$
for certain $b_i\in\F_p$.
Writing $F(X)=\sum_{j=0}^{\infty}a_jX^j$ we deduce that
$ka_k=\sum_{p^i+j=k}b_ia_j$ for all $k$.
In particular, we have $a_1=b_0$, because $a_0=1$.
Replacing $F(X)$ with $F(X)(\bar E_p(a_1X))^{-1}$ we may assume that $a_1=b_0=0$.
The proof will be complete if we show that all $a_j$ with $p\nmid j$ vanish.
Assuming that this is not the case, let $k$ be minimal with $p\nmid k$ and $a_k\neq 0$.
We then have
$0\neq ka_k=\sum_{p^i+j=k}b_ia_j=b_0a_{k-1}=0$, which is the desired contradiction.
Note that this proof works with any field of characteristic $p$ in place of $\F_p$.
\end{rem*}

\begin{proof}[Proof of the Corollary]
We prove that the second condition given in the Corollary is equivalent to the second condition in the Theorem.
It follows from the definition of $E_p(X)$ that
$E_p'(X)=E_p(X)\cdot\sum_{i=0}^{\infty}X^{p^i-1}$.
Hence
$F'(X)=F(X)\cdot\sum_{i=0}^{\infty}X^{p^i-1}$
for every series $F(X)\in\F_p[[X]]$ of the form
$F(X)=\bar E_p(cX)\cdot G(X^p)$
for some $c\in\F_p$ and $G(X)\in 1+X\F_p[[X]]$,
because $G(X^p)$ has null derivative.

To prove the converse, let $F(X)\in\F_p[[X]]$ such that
$F'(X)/F(X)=c\sum_{i=0}^{\infty}X^{p^i-1}$, for some $c\in\F_p$.
Then we have
\[
\frac{d}{dX}\frac{F(X)}{\bar E_p(cX)}=\frac{F'(X)E_p(cX)-F(X)E_p'(cX)}{E_p(cX)^2}=0,
\]
and hence $F(X)/\bar E_p(cX)\in\F_p[[X^p]]$.
\end{proof}

\bibliography{References}

\end{document}